%% file: weyl-group-e8.tex
\documentclass[12pt]{amsart}

\usepackage{amsxtra,amssymb,amsthm,amsmath,amscd,url,xypic,listings}
\usepackage[latin1]{inputenc}
\usepackage{fullpage}
\usepackage{supertabular}

\renewcommand{\leq}{\leqslant}
\renewcommand{\geq}{\geqslant}

\numberwithin{equation}{section}

\newcommand{\Cc}{\mathbf{C}}

\newcommand{\Zz}{\mathbf{Z}}

\newcommand{\Rr}{\mathbf{R}}
\newcommand{\Gg}{\mathbf{G}}

\newcommand{\Qq}{\mathbf{Q}}
\newcommand{\Fp}{\mathbf{F}}
\newcommand{\Tt}{\mathbf{T}}
\newcommand{\T}{\mathbf{T}}

\newcommand{\G}{\mathbf{G}}
\newcommand{\Ee}{\mathbf{E}}

\newcommand{\mods}[1]{\,(\mathrm{mod}\,{#1})}

\newcommand{\ideal}[1]{\mathfrak{{#1}}}

\newcommand{\ra}{\rightarrow}
\newcommand{\lra}{\longrightarrow}

\newcommand{\injecte}{\hookrightarrow}

\newcommand{\fleche}[1]{\stackrel{#1}{\lra}}

\DeclareMathOperator{\Gal}{Gal}

\DeclareMathOperator{\Ad}{Ad}
\DeclareMathOperator{\Lie}{Lie}
\DeclareMathOperator{\ad}{ad}

\newcommand{\eps}{\varepsilon}

\newcommand{\demi}{{\textstyle{\frac{1}{2}}}}

\DeclareMathSymbol{\gena}{\mathord}{letters}{"3C}
\DeclareMathSymbol{\genb}{\mathord}{letters}{"3E}

\theoremstyle{plain}
\newtheorem{theorem}{Theorem}[section]

\newtheorem{lemma}[theorem]{Lemma}

\newtheorem{proposition}[theorem]{Proposition}

\theoremstyle{remark}
\newtheorem{remark}[theorem]{Remark}

\theoremstyle{definition}

\begin{document}

\title{An explicit integral polynomial whose splitting field has
  Galois group $W(\Ee_8)$} 

\author{F. Jouve}
\address{Universit\'e Bordeaux I - IMB\\
351, cours de la Lib\'eration\\
33405 Talence Cedex\\
France}
\email{florent.jouve@math.u-bordeaux1.fr}

\author{E. Kowalski}
\address{Universit\'e Bordeaux I - IMB\\
  351, cours de la Lib\'eration\\
  33405 Talence Cedex\\
  France} \email{emmanuel.kowalski@math.u-bordeaux1.fr}
\thanks{E.K. supported by the A.N.R through the ARITHMATRICS
  project}

\author{D. Zywina} \address{Department of Mathematics, University of
  California, Berkeley, CA 94720-3840, USA}
\email{zywina@math.berkeley.edu}

\dedicatory{Pour les $60$ ans de Henri Cohen}


\begin{abstract}
  Using the principle that characteristic polynomials of matrices
  obtained from elements of a reductive group $\G$ over $\Qq$
  typically have splitting field with Galois group isomorphic to the
  Weyl group of $\G$, we construct an explicit monic integral
  polynomial of degree $240$ whose splitting field has Galois group
  the Weyl group of the exceptional group of type $\Ee_8$.
\end{abstract}

\maketitle

\section{Introduction}

The goal of this paper is to give a concrete explicit example of a
polynomial $P\in\Zz[T]$ such that the Galois group of the splitting
field of $P$ is isomorphic to the group $W(\Ee_8)$, the Weyl group of
the exceptional algebraic group $\Ee_8$. It was motivated by the
construction of such extensions by V\'arilly-Alvarado and
Zywina~\cite{vz} using the Galois action on Mordell-Weil lattices of
some elliptic curves over $\Qq(t)$ which are isomorphic to the root
lattice $\Ee_8$ (this leads in principle to infinitely many such
polynomials, though they are not necessarily easy to write down),
itself based on ideas of Shioda. The existence of such polynomials was
already known from the solution of the inverse Galois problem for Weyl
groups (see the survey of Shioda~\cite{shioda}, or the
paper~\cite{nuzhin} of Nuzhin, as well as~\cite[\S 2.2]{bdeps}
or~\cite[Th. 2]{voskre}).


\begin{theorem}\label{th-1}  
  Let $P\in\Zz[T]$ be the monic polynomial of degree $240$ given by
  $P(T)=T^{120}Q(T+T^{-1})$, where $Q$ is the monic polynomial of
  degree $120$ described by Table~1 in
  Appendix~B.  Then the Galois group of the splitting field of
  $P$ over $\Qq$ is isomorphic to $W(\Ee_8)$.
\end{theorem}

\medskip 
\par
In fact, as we will explain in Proposition~\ref{pr-inf-many}, it is
possible to generalize the construction to obtain infinitely many
(linearly disjoint) examples.  In another direction, although we used
the \textsf{Magma} software~\cite{magma} to construct $P$ (and partly
to prove Theorem~\ref{th-1}), we explain in Appendix~A how it
could be recovered (in principle) ``by hand'', and in particular that
it is quite simple from the point of view of the structure of
reductive algebraic groups.
\par
The basis of the construction is the following principle: if $\G/\Qq$
is a connected reductive algebraic group given as a $\Qq$-subgroup of
$GL(r)$ for some $r\geq 1$, via an injective $\Qq$-homomorphism
$$
\G\fleche{i} GL(r),
$$
and if $g\in \G(\Qq)$ is a ``random'' element, then the Galois group
of the splitting field of $\det(T-i(g))$ (i.e., the characteristic
polynomial of $g$, seen as a matrix through $i$) is typically
isomorphic to the Weyl group $W(\G)$ of $\G$.  Note that such a
principle is in fact pretty close to some of the early methods used
for the study of Lie groups (a ``retour aux sources''), as explained
in the historical notes in~\cite{bourbaki}; in particular, a long time
before the Weyl group was defined in the current manner, \'E. Cartan
(see~\cite{cartan}, in particular pages 50 and following for the case
of $\Ee_8$) determined the Galois group of $\det(T-\ad(X))$ for
a ``general'' $X$ in a simple Lie algebra over $\Cc$ (compare also
with~\cite[\S 8.4, last paragraph]{shioda}, where the same
characteristic polynomial for $\mathfrak{e}_8$ is mentioned and
related to the Mordell-Weil lattices; note those polynomials are not
the same as the ones considered here, e.g, their roots satisfy many
additive relations, whereas ours satisfy multiplicative relations, as
explained in Remark~\ref{rm-relations}).
\par
This principle depends on stating what ``random'' means (and then on
proving the statement!). This was done in~\cite[\S 7]{lsieve} for
elements obtained by random walks
$$
g=\xi_1\cdots \xi_k
$$
in either $SL(r,\Zz)$ ($r\geq 2$, so that $\G=SL(r)$, and $i$ is the
tautological embedding in $GL(r)$) or $Sp(2g,\Zz)$ ($g\geq 1$, so that
$\G=Sp(2g)$ and $i$ corresponds to the standard embedding in
$GL(2g)$): when $k$ is large, the steps of the walk $\xi_j$ being
independently chosen uniformly at random among the elements of a fixed
finite generating set of $\G(\Zz)$, the probability that the splitting
field of $\det(T-g)$ has Galois group different from $W(\Gg)$ is
exponentially small in terms of $k$.
\par
We do not take up the full details of this approach here for the
exceptional group $\Ee_8/\Qq$, though we will come back to this at a
later time in greater generality. What we do is follow the principle
to produce a candidate polynomial. We know that there is an a-priori
embedding of its Galois group in $W(\Ee_8)$, and it turns out (which we
didn't quite expect) that it is possible to check that it is not a
proper subgroup of $W(\Ee_8)$.  

\begin{remark}
  In order to allow easy checking, we have put on the web at the urls
\begin{center}
  \url{www.math.u-bordeaux1.fr/~kowalski/e8pol.gp}
  \url{www.math.u-bordeaux1.fr/~kowalski/e8pol.mgm}
\end{center}
two short files containing definitions of the polynomial above in
\textsf{GP/Pari} and \textsf{Magma}, respectively.  Loading either
will define the variable \verb|pol| to be the polynomial of the
proposition. 
\par
By construction, $P$ is self-reciprocal (so all its roots are
units). Its splitting field turns out to be totally real, and is a
quadratic extension of the splitting field of $Q$. The discriminant of
$P$ is of size about $10^{14952}$, and it is divisible by
\begin{multline*}
  2^{3640}\cdot 3^{300}\cdot 5^{ 30}\cdot 73^{28}\cdot 109^ 2\cdot
  113^ 4\cdot 131^ 4\cdot 331^{ 28}\cdot 419^ { 28}\\
  \cdot 1033 ^ 4\cdot 1103 ^{ 57}\cdot 3307 ^{ 28}\cdot 4649 ^ 4
  \cdot 11467 ^ 4\cdot 629569 ^4\cdot 87087881^4\cdot 508141873^2\\
  \cdot 8321263487^{28} \cdot 58276913161^2 \cdot
  126454995466730813^4\cdot 202992518210175167^{57}
  \\
  \cdot 1644357711723148873333^{28}\cdot
  17520591390337947024593065297057^2,
\end{multline*}
with the cofactor being a square.  Clever use of
\textsf{Pari/GP}~\cite{pari} (as explained by K. Belabas) shows that
the discriminant of the number field of degree $240$ determined by $P$
(i.e., $\Qq[T]/(P)$, not its splitting field) is
$$
1103^{57}\cdot 202992518210175167^{57}
\approx 8.9777\cdot 10^{1159}.
$$
\par
It was also possible to find a polynomial $\tilde{Q}$ such that
$\Qq[T]/(\tilde{Q})\simeq \Qq[T]/(Q)$ with smaller coefficients (by
using the \texttt{polredabs} function), which is available upon
request. 
\end{remark}

\textbf{Notation.} As usual, $|X|$ denotes the cardinality of a
set. For any finite set $R$, $\mathfrak{S}_R$ is the group of all
permutations of $R$, with $\mathfrak{S}_{n}$, $n\geq 1$, being the
case $R=\{1,\ldots, n\}$. We denote by $\Fp_q$ a field with $q$
elements. 
\par
\medskip
\par
\textbf{Acknowledgement.} Many thanks are due to K. Belabas for help
with performing numerical computations (discriminant, basis of the
ring of integers, \texttt{polred}) with the polynomial $P$, etc, and
for explanations of the corresponding functions and algorithms in
\textsf{Pari/GP}; also, thanks to S. Garibaldi for explaining why the
computation with \textsf{GAP} coincides with the one with
\textsf{Magma} (see Appendix~A).

\section{A priori upper bound on the Galois group for $\Ee_8$}
\label{sec-upper}

Let $\Ee_8/\Qq$ be the split group of type $\Ee_8$; it is a simple
algebraic group over $\Qq$ of rank $8$ and dimension $248$. For
information on $\Ee_8$ as a Lie group, we can refer to~\cite{adams}; for
$\Ee_8$ as algebraic group, including proof of existence, abstract
presentation, etc, see, e.g.,~\cite[Ch. 9, Ch. 10, \S
17.5]{springer}. In Appendix~A we also mention a few concrete details.
\par
Contrary to classical groups such as $SL(n)$ or $Sp(2g)$ or orthogonal
groups, which come with an ``obvious'' embedding in a group of
matrices of size comparable with the rank (which is $n-1$ or $g$,
respectively), the smallest faithful representation of $\Ee_8$ is of
dimension $248=\dim \Ee_8$.  More precisely, this is the adjoint
representation
$$
\Ad\,:\, \Ee_8\ra GL(\ideal{e}_8)
$$
where $\ideal{e}_8$ is the Lie algebra of $\Ee_8$, the tangent space at
the identity element with the Lie bracket arising from differentiation
of commutators. This representation is defined over $\Qq$ and given by
$$
g\mapsto T_e(h\mapsto ghg^{-1}),
$$
the differential at the identity element of the conjugation by $g$,
see, e.g,~\cite[I.3.13]{borel}.  The fact that $\Ad$ is injective is
because the center of $\Ee_8$ is trivial (in general, the kernel of the
adjoint representation is the center, in characteristic $0$ at least).
\par
Fix a maximal torus $\T$ of $\Ee_8$ that is defined over $\Qq$ (but not
necessarily split, so that $\T$ is not necessarily isomorphic to
$\G_m^8$ over $\Qq$, but only over some finite extension field; in
fact, the case of interest will be when this field is large).  Let
$X(\T)\simeq \Zz^r$ be the group of characters $\alpha\colon \T \to
\G_m$ (not necessarily defined over $\Qq$).  For each $\alpha\in
X(\T)$, let
\begin{equation*}\label{eq-weight-space}
\ideal{g}_{\alpha}=
\{
X\in \ideal{e}_8\,\mid\,
\Ad(t)\cdot X=\alpha(t)X,\, \text{ for all }t\in \T
\}
\end{equation*}
be the weight space for $\alpha$ in the adjoint representation.  Let
$R(\T,\Ee_8)$ be the set of non-trivial $\alpha\in X(\T)$ with
$\ideal{g}_\alpha \neq 0$; these are called the \emph{roots of $\Ee_8$
  with respect to $\T$}.

\begin{remark}
  It is customary, to view $X(\T)$ as an additive group.  In
  particular, for $\alpha\in R$, the inverse roots $\alpha^{-1}$ is
  denoted $-\alpha$, and $\alpha_1+\alpha_2$ is the root $t\mapsto
  \alpha_1(t)\alpha_2(t)$, etc.
\par
The set $R(\T,\Ee_8)$ is an abstract root system in the space
$V=X(\T)\otimes_\Zz \Rr$; cf.~\cite[Ch.~6]{bourbaki} for definitions.
\end{remark}

The structure theory of reductive groups (see,
e.g.,~\cite[13.18]{borel}) shows that the space $\ideal{g}_\alpha$ is
one dimensional for each root $\alpha\in R(\T,\Ee_8)$, and gives a
direct sum decomposition
\begin{equation} \label{E:lie decomposition}
\ideal{e}_8=\ideal{t}\oplus 
\bigoplus_{\alpha\in R(\T,\Ee_8)}{\ideal{g}_{\alpha}},
\end{equation}
where $\ideal{t}$ is the Lie algebra of $\T$.
\par
From this decomposition, we recover the fact that $|R(\T,\Ee_8)| = \dim
\Ee_8 - \dim \T = 248-8=240$.  The Galois group of $\Qq$ acts naturally
on $X(\T)$: any $\alpha\in X(\T)$ and $\sigma \in
\Gal(\overline{\Qq}/\Qq)$, $\sigma(\alpha)$ is the unique character of
$\T$ such that
$$
\sigma(\alpha(t))= (\sigma(\alpha))(\sigma(t))
$$
for all $t\in\T$. The set of roots $R(\T,\Ee_8)$ is stable under this
action.
\par
Finally, we recall that the \emph{Weyl group of $\Ee_8$ with respect to
  $\T$} is the finite quotient group $W(\T,\Ee_8)=N(\T)/\T$, where
$N(\T)$ is the normalizer of $\T$ in $\Ee_8$.  Since all maximal tori of
a connected linear algebraic group are conjugate (see,
e.g.,~\cite[Th. 6.4.1]{springer}), the Weyl group $W(\T,\Ee_8)$ is
independent of the torus $\T$ up to isomorphism.  We will write
$W(\Ee_8)$ for this abstract group when the choice of torus is
unimportant.
\par
The group $W(\T,\Ee_8)$ acts on the roots by conjugation: for
$w\in N(\T)$, $\alpha\in R(\T,\Ee_8)$, let
\begin{equation}\label{eq-action-w}
(w\cdot \alpha)(t)=\alpha(w^{-1}tw),
\end{equation}
which obviously depends only on the image of $w$ in $W(\T,\Ee_8)$.  This
action is faithful (for instance, because $R(\T,\Ee_8)$ generates the
character group $X(\T)$, and $\T$ is its own centralizer,
see~\cite[13.17]{borel}).
\par
We can now state the main result of this section.
\begin{proposition} \label{P:roots}
  Fix a semisimple element $g\in \Ee_8(\Qq)$, and let $\T$ be any
  maximal torus of $\Ee_8$ that contains $g$.
\par
\emph{(1)} We have the factorization\footnote{\ It is precisely
  because the values of the roots $\alpha$ are the eigenvalues of
  matrices arising from the adjoint representation that the
  terminology \emph{root}, which may seem confusing today, was
  introduced in the historical development of the theory of Lie and
  algebraic groups.}
\begin{equation*}
\det(T-\Ad(g)) = (T-1)^8 \prod_{\alpha \in R(\T,\Ee_8)} (T-\alpha(g)).
\end{equation*}
\par
\emph{(2)} Define the polynomial $P= \det(T-\Ad(g))/(T-1)^8\in
\Qq[T]$, and let $Z\subset \overline{\Qq}$ be the set of roots of $P$.
Assume that $P$ is separable.  Then the map
\begin{equation}\label{E:roots}
\beta\quad
\left\{
\begin{matrix}
 R(\T,\Ee_8)  &\to & Z  \\
\alpha&  \mapsto& \alpha(g)  
\end{matrix}
\right.
\end{equation}
is a bijection which respects the respective
$\Gal(\overline{\Qq}/\Qq)$-actions.
\par
Let $K$ be the splitting field of $P$, i.e., the extension of $\Qq$
generated by $Z$.  Then the Galois action on $R(\T,\Ee_8)$ induces an
injective homomorphism
$$
\phi_g \colon \Gal(K/\Qq) \hookrightarrow W(\T,\Ee_8)
$$
such that for all $\sigma\in \Gal(K/\Qq)$ and $\alpha\in R(\T,\Ee_8)$,
we have
$$
\phi_g(\sigma) \cdot \alpha = \sigma(\alpha).
$$
\end{proposition}

\begin{proof}
  Since $g$ is semisimple, it does lie in a maximal torus $\T$ of $\G$
  (see, e.g.,~\cite[Th. 6.4.5 (ii)]{springer}), and we fix one such
  torus.  The operator $\Ad(g)$ acts as the identity on $\ideal{t}$
  (since conjugation by $g$ is trivial on $\T$) and as multiplication
  by $\alpha(g)$ on each $\ideal{g}_\alpha$, for $\alpha\in
  R(\T,\Ee_8)$.  Therefore from (\ref{E:lie decomposition}), we deduce
  that
\[
\det(T-\Ad(g)) = (T-1)^8 \prod_{\alpha \in R(\T,\Ee_8)} (T-\alpha(g)).
\]
Thus $P$, as defined in the statement of the proposition, is indeed a
polynomial.
\par
Now we assume that $P$ is separable. We first note that
$\alpha(g)\not=1$ for any $\alpha\in R(\T,\Ee_8)$. To see this, we claim
that for any $\alpha$, we can find another root $\alpha'$ such that
$\alpha''=\alpha+\alpha'$ (in additive notation) is also in
$R(T,\Ee_8)$.  Then, since $\alpha'(g)\not=\alpha''(g)$ by assumption,
we obtain $\alpha(g)\not=1$ as desired. From this, in turn, we deduce
(see, e.g.,~\cite[IV.12.2]{borel}) that that $g$ is regular and hence
is contained in a \emph{unique} maximal torus $\T$, which is
necessarily defined over $\Qq$.
\par
Now, to check the claim, one can look at the description of the root
system in Remark~\ref{rm-weyl}, but this is in fact a general property
of any root system $R$ with Dynkin diagram containing no connected
component which is a single point: given $\alpha\in R$, one first
chooses a system $\Delta$ of simple roots such that $\alpha\in\Delta$,
and take $\alpha'$ to be one of the simple roots which are not
perpendicular to $\alpha$ (which exists because of the assumption on
the root system; in other words, $\alpha$ and $\alpha'$ are connected
in the Dynkin diagram of the simple roots; e.g., for $\Ee_8$, if
$\alpha$ corresponds to the vertex labelled $1$ of the Dynkin
diagram~(\ref{eq-dynkin}), one can take $\alpha'$ the root labelled
$3$, etc). Then $(\alpha,\alpha')$ are two simple roots for an
irreducible root system of rank $2$ contained in $R$, and one can
check that $\alpha+\alpha'\in R$ using the classification of those
(see, e.g.,~\cite[9.1.1]{springer}). For $\Ee_8$ (or more generally if
the Dynkin diagram of $R$ has no multiple bond), one can also simply
notice that $s_{\alpha}(\alpha')=\alpha'+\alpha$, where $s_{\alpha}$
is the reflection associated with $\alpha$ (see,
e.g.,~\cite[10.2.2]{springer}).
\par
Coming back to $P$, from the above factorization, we find that the map
$\beta$ is well-defined and surjective, and since
$|R(\T,\Ee_8)|=240=|Z|$, it is therefore bijective.  For each $\alpha\in
R(\T,\Ee_8)$ and $\sigma \in \Gal(\overline{\Qq}/\Qq)$, we have
\[
\beta(\sigma(\alpha)) =\sigma(\alpha)(g) = \sigma(\alpha)(\sigma(g)) =
\sigma( \alpha(g) ),
\]
since $g\in \Ee_8(\Qq)$.  The Galois group $\Gal(K/\Qq)$ acts faithfully
on $Z$ (the permutation action on the roots), so using $\beta$, we
find that $\Gal(K/\Qq)$ acts faithfully on $R(\T,\Ee_8)$, and this
induces an injective group homomorphism
\[
\phi_g \colon \Gal(K/\Qq) \hookrightarrow \mathfrak{S}_{R(\T,\Ee_8)}.
\]
\par
Since $W(\T,\Ee_8)$ acts faithfully on $R(\T,\Ee_8)$, we may naturally
view $W(\T,\Ee_8)$ as a subgroup of $\mathfrak{S}_{R(\T,\Ee_8)}$.  To
conclude, it is thus sufficient to show that the image of $\phi_g$
lies in this subgroup, or in other words, that for every $\sigma \in
\Gal(K/\Qq)$, there exists $w_\sigma \in W(\T, \Ee_8)$ such that
$$
\sigma(\alpha)=w_\sigma\cdot \alpha,\quad\quad
\text{for all $\alpha\in R(\T,\Ee_8)$.}
$$
\par
Fix a split torus $\T_0$ of $\Ee_8$ that is defined over $\Qq$, which
exists since we assumed that our group $\Ee_8$ is split over $\Qq$.
Note that $\T_0$ is split over $K$ and that $\T$ is also.  Indeed, to
check this, it is equivalent to check that the action of
$\Gal(\bar{\Qq}/K)$ on the character group of $\T$ is trivial (see,
e.g.,~\cite[Prop.  13.2.2]{springer}). For this, it suffices to show
that the roots are invariant, since they generate $X(\T)$ (see,
e.g.,~\cite[8.1.11]{springer}, noting that $\Ee_8$ is of adjoint type,
or the description of the root system in Remark~\ref{rm-weyl}). But
for any $\sigma\in \Gal(\bar{\Qq}/K)$, we have
$$
\beta(\sigma(\alpha))=
\sigma(\alpha)(g)=\sigma(\alpha)(\sigma(g))=\alpha(\sigma(g))=
\alpha(g)=\beta(\alpha),
$$
and $\sigma(\alpha)=\alpha$ follows from the injectivity of the map
$\beta$.  
\par
Now the fact that $\T$ and $\T_0$ are both $K$-split implies that
there exists $x\in \Ee_8(K)$ such that $\T=x\T_0x^{-1}$, as proved,
e.g., in~\cite[Th.  15.2.6]{springer}.  Consider then any $\sigma \in
\Gal(K/\Qq)$, and note that $\sigma(x)$ makes sense since $x\in
\Ee_8(K)$. Since both $\T$ and $\T_0$ are defined over $\Qq$, we have
$\T= \sigma(x)\T_0 \sigma(x)^{-1}$ and hence $\sigma(x)x^{-1} \in
N(\T)$.  Let $w_\sigma$ be the element of $W(\T,\Ee_8)$ represented by
$\sigma(x)x^{-1}$. We now claim that $\sigma(\alpha) = w_\sigma \cdot
\alpha$, for all $\alpha \in R(\T,\Ee_8)$, which will finish the proof.
\par
To see this, note that the Galois group $\Gal(K/\Qq)$ acts trivially
on $X(\T_0)$ (because $\T_0$ is split), and that we have an
isomorphism
$$
\gamma
\quad
\left\{
\begin{matrix}
 \T_0 & {\to} &\T\\
 t& \mapsto& xtx^{-1}
\end{matrix}
\right.
$$
which is defined over $K$.  For any $\alpha\in R(\T,\Ee_8)$, we have
 \[
 \sigma(\alpha) = \sigma( \alpha \circ \gamma \circ \gamma^{-1}) = (
 \alpha \circ \gamma) \circ \sigma(\gamma)^{-1},
\]
and then, for all $t\in \T$, we obtain
\begin{equation} \label{E:action on roots} 
\big(\sigma(\alpha)\big)(t)
  = \alpha\big( x \sigma(x)^{-1} t \sigma(x) x^{-1} \big)=\alpha\big(
  \big(\sigma(x)x^{-1}\big)^{-1} t \big(\sigma(x) x^{-1}\big) \big),
\end{equation}
which is the desired conclusion.
\end{proof}

\begin{remark}
  A different approach to Proposition~\ref{P:roots} is sketched (for
  classical groups) in~\cite[App. E]{lsieve}. The one above is more
  direct and intrinsic, and is more amenable to generalizations, but
  we indicate the idea (which can be seen as more down-to-earth):
  given a (regular semisimple) $g\in \Ee_8(\Qq)$, and a fixed
  \emph{split} torus $\T_0$, one considers the set 
$$
X_g=\{t\in \T_0\,\mid\, t\text{ and } g \text{ are conjugate}\}.
$$
\par
This is a non-empty set because $g$ is semisimple, and one shows that
the Weyl group (defined as $N(\T_0)/\T_0$) acts simply transitively by
conjugation on $X_g$; an injection $\Gal(K/\Qq)\ra W(\Ee_8)$ is then
produced by fixing $t_0\in X_g$ and mapping $\sigma$ to $w_{\sigma}$
such that $\sigma(t_0)=w_{\sigma}^{-1}\cdot t_0$. Another small
computation then proves that the permutation of the set of zeros $Z$
obtained from a given $\sigma\in \Gal(K/\Qq)$ is always conjugate to
the permutation of $R(\T_0,\Ee_8)$ induced by $\sigma$.
\end{remark}

\begin{remark}\label{rm-relations}
  Proposition~\ref{P:roots} implies that the zeros of a polynomial
  $\det(T-\Ad(g))$ satisfy many multiplicative relations; indeed, all
  the $240$ zeros are contained in the multiplicative subgroup of
  $\Cc^{\times}$ generated by the $\alpha(t)$ corresponding to eight
  simple roots $\alpha$ (see also~\cite{bdeps} for this type of
  questions, and the next remark if the terminology is unfamiliar).
\end{remark}


\begin{remark}\label{rm-weyl}
  Here are some basic facts on $W(\Ee_8)$ which can be useful to
  orient the reader.
\par
The group $W(\Ee_8)$ is of order $696729600=2^{14}\cdot 3^5\cdot
5^2\cdot 7$, and its simple Jordan-H\"older factors are $\Zz/2\Zz$,
$\Zz/2\Zz$ and the simple group $D_4(\Fp_2)$ (also sometimes denoted
$P\Omega^+_8(2)$, $D_4(2)$, $D_4^+(2)$, or $O_8^+(2)$ as in the Atlas
of Finite Groups~\cite{atlas}), where $D_4$ is the split algebraic
group of type $D_4$ of dimension $28$; this composition series is
essentially already computed by \'E. Cartan in~\cite[p.~50 and
following]{cartan}, working on it as a subgroup of
$\mathfrak{S}_{240}$ (a rather impressive performance). It can be
presented as a Coxeter group (see~\cite[Chapter IV]{bourbaki}) using
eight generators $w_{1}$, \ldots, $w_8$, corresponding to a system of
simple roots $\alpha_1$, \ldots, $\alpha_8 \in R$ (i.e., roots such
that any $\alpha\in R$ can be either represented as integral
combination of the $\alpha_i$ with non-negative coefficient, or its
opposite $\alpha^{-1}$ can be written in this way, but not both),
subject to relations
$$
w_{i}^2=1\,\quad\quad (w_iw_j)^{m(i,j)}=1,\quad\quad 1\leq i<j\leq 8,
$$
where 
$$
m(i,j)=3\text{ if } (i,j)\in \{(1,3),\ (3,4),\ (2,4),\ (4,5),\ (5,6),\
(6,7),\ (7,8)\},
$$
and $m(i,j)=2$ otherwise. (This is encoded in the well-known Dynkin
diagram
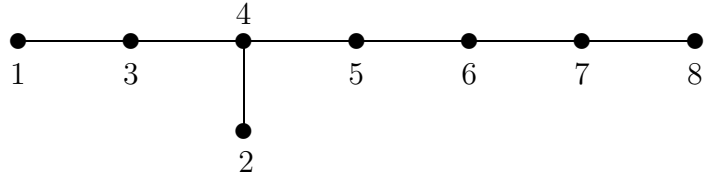
\begin{equation}\label{eq-dynkin}
\unitlength=1mm
\linethickness{0.4pt}
\begin{picture}(142.00,24.00)(25,4)
\put(50.00,22.00){\circle*{2.00}}
\put(65.00,22.00){\circle*{2.00}}
\put(80.00,22.00){\circle*{2.00}}
\put(95.00,22.00){\circle*{2.00}}
\put(110.00,22.00){\circle*{2.00}}
\put(125.00,22.00){\circle*{2.00}}
\put(140.00,22.00){\circle*{2.00}}
\put(80.00,10.00){\circle*{2.00}}
\put(50.00,22.00){\line(1,0){14.00}}
\put(65.00,22.00){\line(1,0){14.00}}
\put(80.00,22.00){\line(1,0){14.00}}
\put(95.00,22.00){\line(1,0){14.00}}
\put(110.00,22.00){\line(1,0){14.00}}
\put(125.00,22.00){\line(1,0){14.00}}
\put(80.00,22.00){\line(0,-1){12.00}}
\put(50.00,19.00){\makebox(0,0)[ct]{$1$}}
\put(65.00,19.00){\makebox(0,0)[ct]{$3$}}
\put(80.00,27.00){\makebox(0,0)[ct]{$4$}}
\put(95.00,19.00){\makebox(0,0)[ct]{$5$}}
\put(110.00,19.00){\makebox(0,0)[ct]{$6$}}
\put(125.00,19.00){\makebox(0,0)[ct]{$7$}}
\put(140.00,19.00){\makebox(0,0)[ct]{$8$}}
\put(80.00,6.00){\makebox(0,0)[c]{\,$2$}}
\end{picture}
\end{equation}
where $m(i,j)=3$ if and only if the vertices labelled $i$ and $j$ are
joined by an edge).
\par
One can also define $W(\Ee_8)$ as the automorphism group of the lattice
$\Gamma_8\subset \Qq^8$ (of rank $8$) generated by $(\demi,\ldots,
\demi)$ and the sublattice
$$
\{(x_1,\ldots, x_8)\in \Zz^8\,\mid\, x_1+\cdots+x_8\equiv 0\mods{2}\},
$$
with the standard bilinear form (see, e.g.,~\cite[V.1.4.3]{serre} for
some more discussion of this lattice, and also~\cite[\S 10]{adams},
where the isomorphism $W(\Ee_8)\simeq \mathrm{Aut}(\Gamma_8)$ is
proved; note many authors studying lattices write $\Ee_8$ for the
lattice instead of the group).  In fact, in the identification of
$W(\Tt,\Ee_8)$, for some maximal torus $\T\subset \Ee_8$, as
$\mathrm{Aut}(\Gamma_8)$, $\Gamma_8$ can be identified with the
character group of $\Tt$, and the roots $R$ are then interpreted as
the $240$ vectors in $\Gamma_8$ with squared-length $2$, namely
\begin{gather*}
  \pm x_i\pm x_j,\quad\quad 1\leq i<j\leq 8,\\
  \demi(\pm x_1\pm x_2\pm \cdots \pm x_8),\quad\quad\text{ with an
    even number of minus signs},
\end{gather*}
the action of $W(\Ee_8)$ on $R$ being the same as the action of the
automorphism group. The lattice $\Gamma_8$ is generated by $R$, with a
basis given for instance by the following eight roots
\begin{gather*}
\demi (x_1+x_2-x_3-x_4-x_5-x_6-x_7-x_8)\\
-x_2+x_3,\quad x_2+x_3,\quad
-x_i+x_{i+1},\text{ for } 3\leq i\leq 7,
\end{gather*}
(which are therefore an example of system of simple roots); see,
e.g.,~\cite[p. 56]{adams}.
\end{remark}

\begin{remark}
See~\cite[\S 7]{shioda} for explicit examples of polynomials whose
splitting fields having Galois groups $W(\Ee_6)$ and $W(\Ee_7)$; they are
much simpler, which can be expected, since $|W(\Ee_6)|=51840=2^7\cdot
3^4\cdot 5$ and $|W(\Ee_7)|=2903040=2^{10}\cdot 3^4\cdot 5\cdot 7$.
Moreover these polynomials have degree $27$, resp.~$28$, which is
smaller than the degrees that would arise from the adjoint
representations, namely $72=78-6$ and $126=133-7$ (this reflect the
fact that there exist faithful representations of the groups $\Ee_6$ and
$\Ee_7$ of simply-connected type in dimension $27$ and $56$).
\end{remark}

\section{Construction of the example}
\label{sec-check}

The polynomial of Theorem~\ref{th-1} is constructed using
\textsf{Magma} (version 2.13-9). We look at the split group $\Ee_8/\Qq$,
and the system of $16$ algebraic generators given by \textsf{Magma},
which come from the Steinberg presentation of reductive algebraic
groups. Precisely (see Appendix~A for some more details and
references), those are the generators $x_{i}=x_{\alpha_i}(1)$, $1\leq
i\leq 8$, of the eight one-parameter unipotent root subgroups
$U_{\alpha_i}$ associated with the simple roots $\alpha_i$ (see,
e.g.,~\cite[8.1.1]{springer}), and the generators
$x_{8+i}=x_{-\alpha_i}(1)$, $1\leq i\leq 8$, of the unipotent
subgroups associated with the negative of the simple roots. The simple
roots are numbered (by \textsf{Magma}) in the usual way described
explicitly, for instance, in~\cite[Ch. VI, \S 4.10]{bourbaki}, and
correspond with the vertices of the Dynkin diagram as in
Remark~\ref{rm-weyl}.
\par
We then construct an element $g$ in $\Ee_8(\Qq)$ by taking the product
of those sixteen generators $x_i$ (in the order above) namely
\begin{equation}\label{eq-expression}
g=x_1\cdots x_{16}=x_{\alpha_1}(1)\cdots x_{\alpha_8}(1)
x_{-\alpha_1}(1)\cdots
x_{-\alpha_8}(1),
\end{equation}
in terms of simple root subgroups; we think of this as a very simple
random walk of length $16$. Then using the adjoint representation of
$\Ee_8$, we compute the matrix $m=\Ad(g)$ (which is in fact in
$SL(248,\Zz)$; in the basis given by \textsf{Magma}, it is a fairly
sparse matrix, with only $6661$ non-zero coefficients among the
$248^2=61540$ entries; the maximal absolute value among the
coefficients is $16$).\footnote{\ Note that we also checked that if we
  construct an element of $\Ee_8(\Qq)$ by taking the product of the $i$
  first generators (in the same order as above) with $1\leq i\leq 15$,
  then the resulting polynomial is not irreducible.}
\par
The characteristic polynomial $\det(T-m)\in\Zz[T]$ is divisible by
$(T-1)^8$ by Proposition~\ref{P:roots}, and the polynomial $P$ of
Theorem~\ref{th-1} is
$$
P=\det(T-m)(T-1)^{-8}.
$$
\par
%
\par
Here are the exact \textsf{Magma} commands to obtain this polynomial
(in a few seconds, this speed depending on fast routines for computing
characteristic polynomials of big integral matrices; neither
\textsf{GAP} nor \textsf{Pari/GP} are able to do this computation
quickly):
\par
\lstset{basicstyle=\small,basicstyle=\ttfamily}
\begin{lstlisting}
A<T>:=PolynomialRing(RationalField()); 
E8:=GroupOfLieType("E8",RationalField()); 
gen:=AlgebraicGenerators(E8); 
rho:=AdjointRepresentation(E8); 
g:=Identity(E8); for i in gen do g:=g*i ; end for; 
m:=rho(g); 
pol:=CharacteristicPolynomial(m) div (T-1)^8; 
\end{lstlisting}
\par
Any decent software package confirms that $P$ is at least irreducible
over $\Qq$ (in particular, its zeros are distinct, as required for the
second part of Proposition~\ref{P:roots}).  Because the roots of $P$
come in inverse pairs, it is possible to write $P=T^{120}Q(T+T^{-1})$
for a unique polynomial $Q\in\Zz[T]$, which we did to shorten a bit
the description of $P$ in Theorem~\ref{th-1}. The irreducibility
of $P$ also implies that $g$ is semisimple: indeed, it suffices to
check that $\Ad(g)$ is diagonalizable, but this is clear because the
minimal polynomial of $\Ad(g)$ has to be $(T-1)P$, and $1$ is not a
zero of $P$.\footnote{\ If $g$ were not semisimple, we could also
  simply argue with its semisimple part, so this is not of great
  importance.}
\par
Now to prove the proposition, let $K$ be the splitting field of $P$, and
$G=\Gal(K/\Qq)$. Thus according to Proposition~\ref{P:roots}, we know
first that $G$ can be identified with a subgroup of $W(\Ee_8)$, and more
importantly that this identification is made in such a way that the
action of $G$ by permutation of the zeros of $P$ in $K$ corresponds to
the action of $W(\Ee_8)$ as a subgroup of $\mathfrak{S}_{240}$ by
permutations of the roots of $\Ee_8$.
\par
This last compatibility is crucial because of the following well-known
fact of algebraic number theory: if $S\in\Zz[T]$ is an irreducible
monic polynomial of degree $d$ with splitting field $F/\Qq$, $H$ the
Galois group of $F/\Qq$ seen as permutation group of the roots of $S$
in $\Cc$, $p$ a prime number such that $S$ factors modulo $p$ in the
form
$$
S\mods{p}=S_1\cdots S_{d},
$$
where $S_i$ is the product of $n_i\geq 0$ distinct monic irreducible
polynomials of degree $i$ in $\Fp_p[T]$, then $H\subset
\mathfrak{S}_d$ contains a permutation with cycle type consisting of
$n_1$ fixed points, $n_2$ disjoint transpositions, etc, and in general
$n_i$ disjoint $i$-cycles.
\par
We apply this to $P$ and $\mathfrak{S}_{240}$, with primes $p=7$ and
$p=11$. We find (again, any decent software package will be able to
factor $P$ modulo $7$ and $11$) that $P\mods{7}$ is the product of $2$
distinct irreducibles of degree $4$, and $29$ distinct irreducibles of
degree $8$, whereas $P\mods{11}$ is the product of $16$ distinct
irreducible polynomials of degree $15$. Hence $G\subset
\mathfrak{S}_{240}$ contains elements of the type
\begin{equation}\label{eq-cycles}
g_{8}=c^{(4)}_1c^{(4)}_2 c^{(8)}_3\cdots c^{(8)}_{31},\quad\quad
g_{15}=d^{(15)}_1\cdots d^{(15)}_{16}
\end{equation}
where the $c^{(\ell)}_i$ (resp. $d_{j}^{(15)}$) are disjoint
$\ell$-cycles (resp. disjoint $15$-cycles). 
\par
In both cases, \textsf{Magma} confirms that such conjugacy classes are
unique in $W(\Ee_8)$ (i.e., there is a single conjugacy class in
$W(\Ee_8)$ with the cycle structure of $g_8$ or $g_{15}$ as permutation
of $R$). 
\par
There are nine conjugacy classes of maximal subgroups in $W(\Ee_8)$,
which are known to \textsf{Magma}. Their indices in $W(\Ee_8)$ are as
follows:
$$
12096,\quad
11200,\quad
2025,\quad
1575,\quad
1120,\quad
960,\quad
135,\quad
120,\quad
2.
$$
\par
Let $M$ be any maximal subgroup; then \textsf{Magma} can also output a
list of the cycle structures, in the permutation action on
$\mathfrak{S}_R\simeq \mathfrak{S}_{240}$, of each conjugacy class of
elements in $M$ (of course, there are sometimes different conjugacy
classes in a given $M$ with the same cycle structure).
\par
Now it turns out, by inspection, that none of the maximal subgroups of
$W(\Ee_8)$ contains elements with the two cycle structures given
in~(\ref{eq-cycles}), and this means that the group $G=\Gal(K/\Qq)$
can not be a subgroup of any of them, and therefore we have
$G=W(\Ee_8)$. 
\par
More precisely, the subgroup of index $2$ is unique and is the kernel
of the restriction of the signature homomorphism $\eps$, which is a
surjective homomorphism
$$
\eps\,:\, W(\Ee_8)\injecte \mathfrak{S}_{240}\ra \Zz/2\Zz,
$$
such that $\eps(g_8)=(-1)^{31}=-1$, $\eps(g_{15})=1$. We see from this
that $G$ is not contained in $\ker\eps$, and hence the only thing to
check to conclude that $G=W(\Ee_8)$ is the fact that none of the maximal
subgroups of index $>2$ contains an element of the class
$g_{15}$. 
\par
This is what we deduced from \textsf{Magma} (but it would be
interesting to have a more conceptual proof; it can also be checked in
the Atlas of Finite Groups~\cite{atlas}, by reducing to the ``big''
simple quotient $D_4(2)=(\ker\eps)/(\text{center})$, for which the
maximal subgroups are listed ``on paper'').
\par
Here are the \textsf{Magma} commands which can be used to construct
$W(\Ee_8)$ and inspect the structure of its maximal subgroups:
\lstset{basicstyle=\small,basicstyle=\ttfamily}
\begin{lstlisting}
W:=WeylGroup(E8); max:=MaximalSubgroups(W);
for m in max do print("----");
  for c in ConjugacyClasses(m`subgroup) do
    print(CycleStructure(c[3]));
  end for;
end for;
\end{lstlisting}
\par
The url \url{www.math.u-bordeaux1.fr/~kowalski/e8check.mgm} contains a
\textsf{Magma} script that lists the maximal subgroups containing
elements of each of the two conjugacy classes (though, as we observed,
checking is only needed for $g_{15}$).

\begin{remark}
  Here are some remarks about this proof, which go in the direction of
  making the objects and arguments more intrinsic and independent of
  an a priori knowledge of the list of maximal subgroups of $W(\Ee_8)$
  (it's not clear if it is reasonable to hope for such a
  proof...). First of all, the conjugacy class of order $15$ is
  particularly symmetric, and we can also prove its uniqueness by pure
  thought. Indeed, it corresponds to the \emph{regular} class of order
  $15$ in $W(\Ee_8)$, as defined by Springer~\cite{springer2}, and
  Springer proved that there is at most one regular conjugacy class of
  a given order in the Weyl group for an irreducible root system
  (see~\cite{springer2}, in particular Theorem 4.1, Proposition 4.10
  and Table 3 in \S 5.4).  Even more precisely, $g_{15}$ is the class
  of the square of the Coxeter elements (e.g.,~\cite[Ch. V, \S
  6]{bourbaki} for the basic properties of the Coxeter element).
\par
Finding the two classes above so easily is somewhat surprising, but it
is not such amazing luck. First, the size of $g_{15}$ is $|W(\Ee_8)|/30$
(again, this can be deduced from Springer's work~\cite[Cor. 4.3,
4.4]{springer2} without invoking any computer check), so by the
Chebotarev density theorem, an extension $L/\Qq$ with Galois group
$W(\Ee_8)$ may be expected to lead to this conjugacy class for roughly
three percent of the primes, which is not negligible.  The class
$g_8$, though less symmetric, is even less surprising from this point
of view: it contains no less than $|W(\Ee_8)|/16$ elements, and is the
largest conjugacy class in $W(\Ee_8)$ (and, as we explained, any odd
conjugacy class would have done just as well for our
argument).\footnote{\ There are $112$ conjugacy classes altogether,
  which are also described explicitly by Carter in~\cite{carter-2}; in
  his notation, $g_8$ is the class with $\Gamma=A_7''$ on p. 56 of
  loc. cit., while $g_{15}$ is the class with $\Gamma=\Ee_8(a_5)$ on p.
  58.}
\par
We state formally the observation we used on subgroups containing
$g_{15}$, as it may prove to be useful for later reference:

\begin{lemma}\label{lm-criterion}
  Let $w_1$, \ldots, $w_8$ be simple reflections generating $W(\Ee_8)$. 
  Let $b=1$ or $2$,
  and let $c=w_1\cdots w_8$ be a Coxeter element in $W(\Ee_8)$. Then any
  proper subgroup of $W(\Ee_8)$ containing an element conjugate to $c^b$
  is contained in the index $2$ subgroup $\ker\eps$.
\end{lemma}

\begin{proof}
  We mentioned that the case $b=2$ is checked unenlighteningly using
  \textsf{Magma}, and then the case $b=1$ follows since a proper
  subgroup containing a conjugate of $c$ contains also a conjugate of
  $c^2$.  (Note that by~\cite[Prop. 4.7]{springer2}, if $b$ is coprime
  with $30$, resp. $15$, then $c^b$ is conjugate to $c$, resp. $c^2$,
  so the lemma holds in fact for any $b$ coprime with $15$.)
\end{proof}
\end{remark}

\section{Infinitely many extensions}

In this section, we show that the construction of the specific
polynomial $P$ also leads easily to infinitely many examples.

\begin{proposition}\label{pr-inf-many}
  Let $\Ee_8/\Zz$ be a model of the split Chevalley group $\Ee_8$ defined
  over $\Zz$, and let $S\subset \Ee_8(\Zz)$ be a symmetric finite
  generating set for $\Ee_8(\Zz)$.  Then
\begin{multline}\label{eq-inf}
\liminf_{k\ra +\infty}
\frac{1}{|S|^k}|\{
(s_1,\ldots, s_k)\in S^k\,\mid\, 
\text{the splitting field of $\det(T-\Ad(s_1\cdots s_k))$}\\
\text{has Galois
  group $W(\Ee_8)$}
\}|>0.
\end{multline}
\par
In particular there are infinitely many $g\in \Ee_8(\Zz)$ for which
$\det(T-\Ad(g))$ has splitting field with Galois group isomorphic to
$W(\Ee_8)$.
\end{proposition}

As explained before, one can expect a much stronger result (the
left-hand side of~(\ref{eq-inf}) should be $\geq 1-C\exp(-ck)$ for
some $c>0$, $C\geq 0$), but checking this would involve a deeper
analysis of the finite groups $\Ee_8(\Fp_q)$, which we defer to another
time. Also, it should be possible to prove in this manner the
existence of infinitely many polynomials with (globally) linearly
disjoint splitting fields with Galois group $W(\Ee_8)$ (this is already
known, see~\cite[Th. 7.1]{shioda}).

\begin{proof}
  First of all, the fact that $\Ee_8(\Zz)$ is finitely generated (hence
  $S$ exists) is a standard property of Chevalley
  groups. 
\par
Let $g$ be the element of $\Ee_8$ in the proof of
Theorem~\ref{th-1}; it turns out that $g\in \Ee_8(\Zz)$ (this is
clear from~(\ref{eq-expression}) and the fact that \textsf{Magma}
constructs a group defined over $\Zz$). Let
$P=\det(T-\Ad(g))(T-1)^{-8}$.
Now, we claim that for any $h\in \Ee_8(\Zz)$, if $h$ is conjugate to $g$
modulo $p$ for $p=7$ and $p=11$
(where congruences refer to the reduction maps $\Ee_8(\Zz)\ra
\Ee_8(\Fp_p)$, or to congruences of matrices after applying $\Ad$, and
conjugation is in $\Ee_8(\Fp_p)$), then the Galois group of the
splitting field of $Q=\det(T-\Ad(h))(T-1)^{-8}$ must be $W(\Ee_8)$.
\par
Indeed, let $h_s\in \Ee_8(\Qq)$ be the semisimple part of $h$ (see,
e.g.,~\cite[I.4.4]{borel}); we also have
$$
Q=\det(T-\Ad(h_s))(T-1)^{-8}.
$$
\par
For $p=7$ and $p=11$, we have $Q\equiv P\mods{p}$, and since $P$ has
distinct roots modulo $11$, not including $1\in \Fp_{11}$ (it has only
irreducible factors of degree $15$), these conditions imply that $h_s$
must be regular semisimple, and that $Q$ has distinct roots.
\par
Finally, the Galois group of the splitting field of $Q$ will contain
elements of the same conjugacy classes $g_8$ and $g_{15}$ discussed in
the proof of Theorem~\ref{th-1}, and hence by
Proposition~\ref{P:roots}, it will have to be isomorphic to $W(\Ee_8)$.
\par
Now let
$$
H=\Ee_8(\Fp_7)\times \Ee_8(\Fp_{11}).
$$
\par
Because the $\Ee_8(\Fp_p)$ are distinct non-abelian simple groups for
all $p\geq 2$ (this is due to Chevalley~\cite{chevalley}), the
reduction map $\Ee_8(\Zz)\fleche{\pi} H$ is surjective.  Indeed, the
individual reduction maps $\Ee_8(\Zz)\ra \Ee_8(\Fp_p)$ are onto,
because the algebraic generators $x_{\alpha}(1)$ in $\Ee_8(\Zz)$
associated with the roots of $\Ee_8$ (with respect to a split maximal
torus) reduce to the corresponding generators of $\Ee_8(\Fp_p)$ (see,
e.g.,~\cite[\S 6]{steinberg} for the fact that the elements
$x_{\alpha}(1)$ generate the group of rational points of a simple
split Chevalley group over a prime field; this can also be checked for
$p=7$, $11$ with \textsf{Magma}'s \texttt{Generators()} command), and
one can apply the classical Goursat lemma to the image of $\pi$ (a
proper subgroup of $G_1\times G_2$, where $G_i$ are non-abelian simple
groups, which surjects to $G_1$ and $G_2$, is the graph of an
isomorphism $G_1\ra G_2$).
\par
Then it is a standard fact about random walks on finite groups
(``convergence to the invariant distribution of reversible, aperiodic,
irreducible, finite Markov chains'') that we have
$$
\lim_{k\ra +\infty}
\frac{1}{|S|^k}|\{
(s_1,\ldots, s_k)\in S^k\,\mid\, 
\pi(s_1\cdots s_k)\text{ is conjugate to } 
(g,g)\in H\}|=\frac{|C|}{|H|},
$$
where $C\subset H$ is the conjugacy class of $(g,g)$ (see the
discussion in~\cite[Th. 2.1, \S 2.2]{saloff-coste} and~\cite[Chapter
7]{lsieve}; in our case, the aperiodicity follows from the symmetry of
$S$, and the fact that there is no non-trivial homomorphism
$\Ee_8(\Fp_p)\ra \Zz/2\Zz$).
\par
It follows from the two observations above that the proposition holds
with the precision that the liminf is $\geq |C||H|^{-1}$ (which,
however, is very small, roughly $10^{-15}$). Finally, although
distinct ``words'' $(s_1,\ldots,s_k)$ may lead to the same element,
the result clearly implies the existence of infinitely many distinct
$h$ with the desired property (e.g., because if there were only a
finite list $(h_1,\ldots, h_N)$ of such, we could repeat the argument
with additional congruences $s_1\cdots s_k\not\equiv h_i\mods{p_1}$
where $p_1$ is a prime chosen so that the $h_i$ do not represent all
classes modulo $p_1$, e.g., $p_1>N$, to obtain a contradiction).
\end{proof}

\section*{Appendix A: intrinsic characterization of the polynomial}

We now build on~(\ref{eq-expression}) to explain in detail how the
definition (and computation) of $P$ may be phrased in such a way that
it does not depend on any choice or implementation detail in
\textsf{Magma}'s code (which may, in particular, vary from version to
version). So, in principle, it would be possible to compute $P$ by
hand using only printed references (such as~\cite{steinberg}
or~\cite{cmt}). More practically, other programs can be used to check
the computation.
\par
To make things clearer, we denote here by $\Ee_8^m/\Qq$ the split group
of type $\Ee_8$ over $\Qq$ given by \textsf{Magma}. Associated with it
are a maximal torus $\Tt^m\subset \Ee_8^m$, split over $\Qq$, the set of
roots $R$ associated with $\Tt^m$, and a certain choice $\Delta\subset
R$ of simple roots. Those are enumerated
$$
\Delta=\{\alpha_1,\ldots,\alpha_8\}
$$
as dictated by the Dynkin diagram: the roots $\alpha_i$ and $\alpha_j$
are not orthogonal, with respect to a $W(\Ee_8)$-invariant inner product
on $X(\Tt^m)\otimes\Rr$, if and only if the vertices $i$ and $j$ of
the Dynkin diagram are connected.
\par
For each $\alpha\in R$, there is a one-parameter unipotent root
subgroup $U_{\alpha}$ which is the image of a non-trivial homomorphism
$$
x_{\alpha}\,:\, \Gg_a\ra \Ee_8^m
$$
which is defined over $\Qq$ and such that
$$
t^{-1} x_{\alpha} (u) t=x_{\alpha}(\alpha(t)u)
$$
for $t\in \Tt^m$ and $u\in\Gg_a$. The generators giving $g\in
\Ee_8^m(\Qq)$ in~(\ref{eq-expression}) are $x_i=x_{\alpha_i}(1)$ for
$1\leq i\leq 8$ and $x_{8+i}=x_{-\alpha_i}(1)$ for $1\leq i\leq 8$.
\par
To compute $\Ad(g)$, since $\Ad$ is an homomorphism, one needs to
compute $\Ad(x_{\alpha}(u))$ for $\alpha\in R$ and $u\in \Qq$. Now we
have an induced map between Lie algebras
$$
\Lie(\Gg_a)\fleche{dx_{\alpha}} \Lie(\Ee_8^m).
$$
\par
Define $e_{\alpha}=dx_{\alpha}(1)\in \Lie(\Ee_8^m)$; this is a generator
of the root space $\mathfrak{g}_{\alpha}$ associated to
$\alpha$. Because the image of $x_{\alpha}$ is unipotent,
$\ad(e_{\alpha})$ is a nilpotent endomorphism of $\Lie(\Ee_8^m)$, where
$\ad$ is the adjoint representation at the Lie algebra level (so that
$\ad(X)$ maps $Y$ to $[X,Y]$, where $[X,Y]$ is the Lie bracket). Then
we have the formula
\begin{equation}\label{eq-exp}
\Ad(x_{\alpha}(u))=\exp(u\ad(e_{\alpha})),
\end{equation}
where the exponential, which can be interpreted by the usual power
series as an exponential of matrix, is in fact a polynomial in
$u\ad(e_{\alpha})$ since $\ad(e_{\alpha})$ is nilpotent
(see~(\ref{eq-nilp}) below). (This can be proved purely algebraically,
but we may also extend scalars to $\Rr$, and see that both sides
represent smooth functions of $u\in \Rr$ into
$GL(\Lie(\Ee_8^m)\otimes\Rr)$ which satisfy the same ordinary
differential equation $\tfrac{dy}{du}=\ad(e_{\alpha})y$ and which take
the same value at $0$).
\par
Thus to compute $\Ad(g)$, it is enough to compute the endomorphisms
$\ad(e_{\alpha})$ for $\alpha\in R$. But since a basis of the Lie
algebra is made of a basis (say $h_1,\ldots, h_8$) of the Lie algebra
of the torus $\Tt^m$, and the $e_{\alpha}$ for $\alpha\in R$, this
amounts in turn to being able to compute the brackets
$[e_{\alpha},e_{\beta}]$ for $\alpha$, $\beta\in R$ and
$[e_{\alpha},h_i]$ for all $\alpha$ and $i$.
\par
It turns out that those brackets are explicitly known and depend only
on the ``abstract'' root system $R$ except for
\begin{gather*} 
  [e_{\alpha},e_{\beta}]=c(\alpha,\beta)e_{\alpha+\beta}
\end{gather*}
where $c(\alpha,\beta)\in \Qq^{\times}$. Those $c(\alpha,\beta)$ are
known as the structure constants for the Lie algebra; in fact, when
the group comes from a group scheme defined over $\Zz$ (as is the case
of $\Ee_8^m$), we have $c(\alpha,\beta)\in \{\pm 1\}$. At the level of
the group, the structure constants occur in the commutator relations
$$
[x_{\alpha}(u),x_{\beta}(v)]=x_{\alpha+\beta}(c(\alpha,\beta)uv)
$$
for $\alpha$, $\beta\in R$ with $\alpha+\beta\in R$ and $u$,
$v\in\Gg_a$ (the simple form of this relation is due to the fact that
the root system of $\Ee_8$ is an example of \emph{simply laced} root
system, see, e.g.,~\cite[\S 10.2]{springer}).
\par
Note in passing that the other brackets imply in particular that
$\ad(e_{\alpha})$ is nilpotent of order $2$, so that~(\ref{eq-exp})
becomes
\begin{equation}\label{eq-nilp}
\Ad(x_{\alpha}(u))=\exp(u\ad(e_{\alpha}))=\mathrm{Id}+u\ad(e_{\alpha})+
\frac{u^2}{2}\ad(e_{\alpha})^2,
\end{equation}
see, e.g.,~\cite[10.2.7]{springer} or~\cite[\S 3]{cmt}.
\par
So the endomorphism $\Ad(g)$ is easily computable from the knowledge
of the structure constants.  However, matters are somewhat complicated
from then on by the fact that there is no absolutely canonical choice
of the $c(\alpha,\beta)$. Still, as described for instance in~\cite[\S
2.3, \S 3]{cmt}, once a certain total order has been put on the root
system, there exists a certain set of \emph{extraspecial pairs}
$(\alpha,\beta)$, precisely $112$ of them, for which $c(\alpha,\beta)$
can be chosen arbitrarily in $\{\pm 1\}$, and then all other structure
constants are uniquely determined. 
\par
Thus to describe unambiguously our endomorphism $\Ad(g)$, it suffices
to describe the extraspecial structure constants in
$\Lie(\Ee_8^m)$. \emph{These are defined to all be $+1$}. This can be
checked by the following \textsf{Magma} commands:
\par
\lstset{basicstyle=\small,basicstyle=\ttfamily}
\begin{lstlisting}
L:=LieAlgebra(GroupOfLieType("E8",RationalField()));
ExtraspecialSigns(RootDatum(L));
\end{lstlisting}
\par
This already provides a way to construct from scratch, in principle,
the polynomial $P$ of Theorem~\ref{th-1}. However, there is an even
stronger ``unicity'' feature, which was explained to us by Skip
Garibaldi: for any choice of generators $x_{\alpha}(1)$ of the
unipotent root subgroups (of a split group $\Ee_8$ of type $E_8$ over
$\Zz$, with split maximal torus $\Tt/\Zz$ and simple roots
$\Delta=\{\alpha_1,\ldots,\alpha_8\}$), determining the generators
$x_i$ as above, the element
$$
g=x_1\cdots x_8x_9\cdots x_{16}
$$
has \emph{the same} characteristic polynomial. The point is that the
elements $x_i$ are determined up to sign from the choice of the simple
roots, hence the possible changes are determined by a vector
$\boldsymbol{\eps}=(\eps_1,\ldots,\eps_8)\in \{\pm 1\}^8$ of signs,
and the possible elements $g$ that can be obtained are, relative to a
fixed group $\Ee_8/\Zz$, of the form
$$
g_{\boldsymbol{\eps}}=x_{\alpha_1}(\eps_1)\cdots
x_{\alpha_8}(\eps_8)
x_{-\alpha_1}(\eps_1)\cdots
x_{-\alpha_8}(\eps_8)
$$
\par
Now it turns out that there exists an element
$t_{\boldsymbol{\eps}}\in\Tt$, depending only on those signs, such
that
$$
x_{\alpha_i}(\eps_i)=t_{\boldsymbol{\eps}}x_{\alpha_i}(1)t_{\boldsymbol{\eps}}^{-1}
$$
for all simple roots $\alpha_i$ (this follows, e.g.,
from~\cite[VIII.5.2, Cor. 3]{bourbaki2}). Then a simple computation
(which can be done in $SL(2)$, because it only concerns a root and its
negative) shows that we also have
$$
x_{-\alpha_i}(\eps_i)=t_{\boldsymbol{\eps}}x_{\alpha_i}(1)t_{\boldsymbol{\eps}}^{-1}
$$
and therefore we also have
$$
g_{\boldsymbol{\eps}}=t_{\boldsymbol{\eps}}gt_{\boldsymbol{\eps}}^{-1},
$$
so that all $\Ad(g_{\boldsymbol{\eps}})$ are conjugate and have the
same characteristic polynomial.
\par
We implemented this strategy using the \textsf{GAP} system~\cite{gap},
version $4.4.9$, which knows about Lie algebras (but not algebraic
groups), and has different structure constants than those of
\textsf{Magma} (for instance, there is an extraspecial pair
  $(\alpha_1,\alpha_3)$, and
  $[e_{\alpha_1},e_{\alpha_3}]=e_{\alpha_{1}+\alpha_3}$ for
  \textsf{Magma}, while
  $[e_{\alpha_1},e_{\alpha_3}]=-e_{\alpha_{1}+\alpha_3}$ for
  \textsf{GAP}).
The recipe above, as it should, leads to a matrix with the same
polynomial $P$ as in Theorem~\ref{th-1} (note that the
\texttt{CharacteristicPolynomial} function in \textsf{GAP} is not up
to the task of computing $P$ from the matrix in a reasonable amount of
time, so we did this last check using \textsf{Magma} again, though one
could also check modulo sufficiently many small primes to ensure the
result by the Chinese Remainder Theorem).  Here are the commands to
produce this matrix: 
\lstset{basicstyle=\small,basicstyle=\ttfamily}
\begin{lstlisting}
L:=SimpleLieAlgebra("E",8,Rationals);
d:=[];
for i in [1..248] do 
  Append(d,[AdjointMatrix(Basis(L),Basis(L)[i])]); 
od;
a:=[];  
am:=[];
for i in [1..8] do
  a[i]:=IdentityMat(248)+d[i]+1/2*d[i]*d[i]; 
  am[i]:=IdentityMat(248)+d[120+i]+1/2*d[120+i]*d[120+i]; 
od;
m:=a[1]*a[2]*a[3]*a[4]*a[5]*a[6]*a[7]*a[8];
m:=m*am[1]*am[2]*am[3]*am[4]*am[5]*am[6]*am[7]*am[8]; 
\end{lstlisting}
\par
Since \textsf{GAP} is Open Source, this computation can (or could) be
checked in complete detail, guaranteeing the correctness of
Theorem~\ref{th-1}.

\section*{Appendix B: coefficient table}

We conclude with an Appendix listing the table of coefficients of the
polynomial $Q$ such that $P(T)=T^{120}Q(T+T^{-1})$.
\par
\input e8polbis.tex
\par

\end{document}

%% file: e8polbis.tex
\begin{center}
\tablehead{
\emph{Degree} $i$ &\emph{Coefficient of } $T^i$ \\
\hline
}
\begin{supertabular}{c|c}
$0$  &   $\scriptstyle{365587894983967922854560421106658794835269162025801581455381270108994772086354097689969}$ \\
$1$  &   $\scriptstyle{-11188764671313743052104852260462353867603756780241598167388311775144605678410262867332448}$ \\
$2$  &   $\scriptstyle{169253696238399029559192135798369681596869020592118579039666548219126339368278708130365896}$ \\
$3$  &   $\scriptstyle{-1687159920262524494571824891028833720039695470637056956223462293372163529290126431304585100}$ \\
$4$  &   $\scriptstyle{12466538870569350428117512482674638738192598391028225510113876433671029848121166234252784168}$ \\
$5$  &   $\scriptstyle{-72826947156697363455035723423426971890324922381314772703337692580409097282009564861222315768}$ \\
$6$  &   $\scriptstyle{350331529672673601609711561533019721386090515420775351136161318553244630344247566318282845504}$ \\
$7$  &   $\scriptstyle{-1427238565615616836874966027280849892931467391947083350631098008579547451949365512165445324548}$ \\
$8$  &   $\scriptstyle{5026354438587350042393019188859322619309223067912455590765316222074769006654715395492167534064}$ \\
$9$  &   $\scriptstyle{-15543242692747890030945651958468728793028138424448388052344884541152421476286890160281606087912}$ \\
$10$  &   $\scriptstyle{42727722856586577963125858580310302599252714404522901925297264385654281164703726634073363638112}$ \\
$11$  &   $\scriptstyle{-105456673870095173711703163143669727224697814378273599459423442778596360012854640538260794921080}$ \\
$12$  &   $\scriptstyle{235607817095755199368858190149899878856163978515561431455451403131557914807361950016856603119900}$ \\
$13$  &   $\scriptstyle{-479770149264907609591636867066530148444729317862800608632937138608427032184610743855219428626708}$ \\
$14$  &   $\scriptstyle{895635623173061803600226851312211772312504932418933306086383616156548261949688014524531333153508}$ \\
$15$  &   $\scriptstyle{-1540472555790089754392061652689278932257505758853194366507839523329666453963515188071478275995640}$ \\
$16$  &   $\scriptstyle{2451788710586871745589664228861539890944914421283328064352160762645878635406460765371848464480340}$ \\
$17$  &   $\scriptstyle{-3624640505901638188009697337857954483473385109589930872949558184551987911180781360042790583083604}$ \\
$18$  &   $\scriptstyle{4994024072126951505126399098824620904652453895567294295762673834554854694248727868906577016537540}$ \\
$19$  &   $\scriptstyle{-6431721637672973764253559554028129752608859470080682240640404627541331731173957688036200659035092}$ \\
$20$  &   $\scriptstyle{7763241468623979545008124501105553771540039037756350364261700519310257100760280299151055542732492}$ \\
$21$  &   $\scriptstyle{-8802942880381499388245651070259273872571338870198906824716298766585276914830555824919832373777496}$ \\
$22$  &   $\scriptstyle{9397428905594516743427299464872348062738848222907978714295289389671821494398632421489330984834096}$ \\
$23$  &   $\scriptstyle{-9462921861838231504232648180683566209060687542850222756830572873896415321834823788805352470036944}$ \\
$24$  &   $\scriptstyle{9004042910145318452280462187844950546979429562212010153159452511976316511237893389182144564071816}$ \\
$25$  &   $\scriptstyle{-8108466777447973350625143866020889145905001781484033233303172750662906603502257635040763697120800}$ \\
$26$  &   $\scriptstyle{6920877993691644199503325490868076639226322468228894863452356759758353696538512384603684679212236}$ \\
$27$  &   $\scriptstyle{-5606393286650355050632445312392545648077485412935086399427015013222455115423342550668229228526380}$ \\
$28$  &   $\scriptstyle{4315552971543918712292689802749919612866096319783016002162626109826378082792279094813015908629560}$ \\
$29$  &   $\scriptstyle{-3160147292945573961625584424097717920370225757650282483390682408815821565443532117419343460827840}$ \\
$30$  &   $\scriptstyle{2203660581121260808668034624515265900353324373444058449372337614106044024846226145678411752103944}$ \\
$31$  &   $\scriptstyle{-1464744343303019546273900066171788739473025422537466727740192615776481945774661361874676553542888}$ \\
$32$  &   $\scriptstyle{928837747351900139847510911539219291666880836492600960301761612482212279686785348971171633827678}$ \\
$33$  &   $\scriptstyle{-562380987452023460657841067338702972535374341028182139165378619782789496619429231331416465467136}$ \\
$34$  &   $\scriptstyle{325356235824906304914262446976176796571764856456686279807557436317349277524804415206207943829668}$ \\
$35$  &   $\scriptstyle{-179980581397996318804833400684252642822685183748168637880480218232713417673333008007141832612656}$ \\
$36$  &   $\scriptstyle{95259358673342811962577374968024344116759528191118197330901994055953518185301961938058278687152}$ \\
$37$  &   $\scriptstyle{-48268244673830087725415910980532812904398769665363772047167723957702501678922709545345961314952}$ \\
$38$  &   $\scriptstyle{23427392442378901750354958702292727966207504749183723357421289869802194332649685331485983388676}$ \\
$39$  &   $\scriptstyle{-10897163294836081802023740205919221704582045121673085440709523835797116373484015837949721485760}$ \\
$40$  &   $\scriptstyle{4859945481426652393657932982123404648051819351769400595806382044197187788576701075955522665148}$ \\
$41$  &   $\scriptstyle{-2079041289572381880981851137517243250745330223456680067657643451728129396120093304826885539344}$ \\
$42$  &   $\scriptstyle{853451569340348395409517113261975739247889854460920295141494937564435136719739663763259170648}$ \\
$43$  &   $\scriptstyle{-336306101614439964957860525247000005857956841159545332435837005318322034194039653447757128580}$ \\
$44$  &   $\scriptstyle{127254544703467199806219784585607568957868365018715949880037719179469076397221511835490064520}$ \\
$45$  &   $\scriptstyle{-46251344543607872597100973940088670916183046876614132396624422211322670222925432802603348764}$ \\
$46$  &   $\scriptstyle{16151229449214495754070501316928248080301470592125995573153793317058666430416170267465079692}$ \\
$47$  &   $\scriptstyle{-5420296230053687740961180866132183229319212804267569609857245172190699563010303045911368944}$ \\
$48$  &   $\scriptstyle{1748523996015819207122503672053973801759807116329056619663716454198769892066785902652554878}$ \\
$49$  &   $\scriptstyle{-542294947655157376138526043814456732351589519118047194151659330168304308215153590486208844}$ \\
$50$  &   $\scriptstyle{161729358620880031971279266977475707554378016762290238000893567282554531020744559074376324}$ \\
$51$  &   $\scriptstyle{-46386918898623721747469881003020417214566850265053511450799223706659896381395790959743804}$ \\
$52$  &   $\scriptstyle{12797051499162896398440721779023796724984325295916014796039628931842465948343894845111380}$ \\
$53$  &   $\scriptstyle{-3396073129112728218339337536714697724329251745742194506364950008794369800942410905152136}$ \\
$54$  &   $\scriptstyle{867027209473064571745438316617176614625159965185956598478937806132721701397948869458900}$ \\
$55$  &   $\scriptstyle{-212962875904766156808053958032907333576950046165625774677667251904787655643757505176324}$ \\
$56$  &   $\scriptstyle{50327885532415893880483155785349071253807731172747773207451521031252671732632230388978}$ \\
$57$  &   $\scriptstyle{-11443438667187908543351570986876579176110747490338119245195350888150296055016975582004}$ \\
$58$  &   $\scriptstyle{2503504763126665005017589126659165629183801291281925243578712171462099771849619254744}$ \\
$59$  &   $\scriptstyle{-526958944944907268493639103677877316954184364418982668079438920193188411899643568508}$ \\
$60$  &   $\scriptstyle{106715051592438328930689239004654253545573721719451360040168657740453189940915192314}$ \\
$61$  &   $\scriptstyle{-20790779910335793865471457635357833163154112951516080880730395764735512013025652796}$ \\
$62$  &   $\scriptstyle{3896541990320087145594716450422460409978287599082283693485952309365925318005644500}$ \\
$63$  &   $\scriptstyle{-702441672053715648053554868151235620650214741238958432207947706491608211081312400}$ \\
$64$  &   $\scriptstyle{121790629393362604544357589123597436252606558743434815262003043923897464121167197}$ \\
$65$  &   $\scriptstyle{-20306355805901050837901846972656773418756565045468863673412540172210337271902340}$ \\
$66$  &   $\scriptstyle{3255350507442216588937835814600680908204323769283956961525023238451439544941760}$ \\
$67$  &   $\scriptstyle{-501690059905683053720524741882157670117822132733245580070161190249583644847028}$ \\
$68$  &   $\scriptstyle{74312275258414477271317122327593160444507387370961124773453087824921327894416}$ \\
$69$  &   $\scriptstyle{-10577404455255469741143362965021672227782124046356528509651977109935992639428}$ \\
$70$  &   $\scriptstyle{1446398328284711960698545660047973219044522333438207563556095761539846607256}$ \\
$71$  &   $\scriptstyle{-189965195590509120305163972290005629377791264599483176752468398233763249484}$ \\
$72$  &   $\scriptstyle{23956016062333717865546677021452562824319209365764369018142153751707952028}$ \\
$73$  &   $\scriptstyle{-2899855516264429489706362148867135560476890463566496024776316822467683664}$ \\
$74$  &   $\scriptstyle{336833599313788500559679181776716047210388540402180437654119288339393740}$ \\
$75$  &   $\scriptstyle{-37529836525283747105371855723169247756700129909528157240024250001998788}$ \\
$76$  &   $\scriptstyle{4009530152235809687640827672341358465891595297597478301291898189654876}$ \\
$77$  &   $\scriptstyle{-410569812715360439015853397370594547170545487427635936414500045030652}$ \\
$78$  &   $\scriptstyle{40277795973917435115762646297317343617205272567786230485171910510652}$ \\
$79$  &   $\scriptstyle{-3783768552039089682586618639655876942103322356223597420264279396724}$ \\
$80$  &   $\scriptstyle{340207987956327503878827250457263328265317830300159719894983206382}$ \\
$81$  &   $\scriptstyle{-29261155514270620208017844865051944340304021710386304703117145200}$ \\
$82$  &   $\scriptstyle{2406110405248787166244163214026562342878677865734412916133200336}$ \\
$83$  &   $\scriptstyle{-189038373443757680942334919469355568250711998730640952299841492}$ \\
$84$  &   $\scriptstyle{14181101939087673375192221173276978973988251595771099469911248}$ \\
$85$  &   $\scriptstyle{-1015059691711179211122508131842675818892409168967000131556504}$ \\
$86$  &   $\scriptstyle{69273982784497519473037155287005412956159996851830481504940}$ \\
$87$  &   $\scriptstyle{-4504005163576022050795556544568918073128228705604683414856}$ \\
$88$  &   $\scriptstyle{278744922765563512122592853176105123936966933933756780006}$ \\
$89$  &   $\scriptstyle{-16405810574710923918958669050194823065267996065080494420}$ \\
$90$  &   $\scriptstyle{917373448095043315139701983319126630369060417226308240}$ \\
$91$  &   $\scriptstyle{-48685180731396433963474339952113579549797184903456932}$ \\
$92$  &   $\scriptstyle{2449403736861952764502194313954481630464621913377362}$ \\
$93$  &   $\scriptstyle{-116683889096980060401747351223521809696783299718096}$ \\
$94$  &   $\scriptstyle{5256318681823135646376194757383743386905542401984}$ \\
$95$  &   $\scriptstyle{-223594517259941108513934941256272336603486583080}$ \\
$96$  &   $\scriptstyle{8967837569963007481063656618084357384523513289}$ \\
$97$  &   $\scriptstyle{-338566611781299224336234061335519673555729968}$ \\
$98$  &   $\scriptstyle{12010206026853145744238047620715510383030860}$ \\
$99$  &   $\scriptstyle{-399536336273091783248772485529679824619372}$ \\
$100$  &   $\scriptstyle{12437443190985072692004053616601920326304}$ \\
$101$  &   $\scriptstyle{-361454027338752388080364343617032962704}$ \\
$102$  &   $\scriptstyle{9781270532352502601151919455054150468}$ \\
$103$  &   $\scriptstyle{-245758896101673421544480044149574716}$ \\
$104$  &   $\scriptstyle{5714846913247894902773082359642858}$ \\
$105$  &   $\scriptstyle{-122552946275635592994659815458660}$ \\
$106$  &   $\scriptstyle{2413849634493902632445738578404}$ \\
$107$  &   $\scriptstyle{-43467840407415668458306853984}$ \\
$108$  &   $\scriptstyle{711884504814065975117065754}$ \\
$109$  &   $\scriptstyle{-10538669572997700747046736}$ \\
$110$  &   $\scriptstyle{140021126816597308605612}$ \\
$111$  &   $\scriptstyle{-1655565532193307303324}$ \\
$112$  &   $\scriptstyle{17242140511966984109}$ \\
$113$  &   $\scriptstyle{-156184748605164508}$ \\
$114$  &   $\scriptstyle{1211012431626440}$ \\
$115$  &   $\scriptstyle{-7871527038772}$ \\
$116$  &   $\scriptstyle{41688975082}$ \\
$117$  &   $\scriptstyle{-172657460}$ \\
$118$  &   $\scriptstyle{524076}$ \\
$119$  &   $\scriptstyle{-1036}$ \\
$120$  &   $\scriptstyle{1}$ \\
\end{supertabular}
\end{center}